\documentclass[aps,reprint,nofootinbib]{revtex4-1}
\usepackage{graphicx} \usepackage{amsfonts} 
\usepackage{times} \usepackage{subfigure}
\usepackage{amssymb}\usepackage{amsmath} 
\usepackage[usenames,dvipsnames]{color}

\newcommand{\beq}{\begin{equation}}     \newcommand{\eeq}{\end{equation}}
\newcommand{\beqa}{\begin{eqnarray}}    \newcommand{\eeqa}{\end{eqnarray}}
\newcommand{\bde}{\begin{description}}  \newcommand{\ede}{\end{description}}
\newcommand{\ben}{\begin{enumerate}}    \newcommand{\een}{\end{enumerate}}

\newcommand{\R}{{\mathbb{R}}}

\newcommand{\eqn}[1]{\beq{ #1 }\eeq}

\newcommand{\inv}[1]{{\frac{1}{#1}}}

\newcommand{\inRbracket}[1]{{\left({#1}\right)}}
\newcommand{\inSbracket}[1]{{\left[{#1}\right]}}

\usepackage{graphicx, rotating}
\usepackage{epsfig}
\usepackage[abs]{overpic}
\usepackage{psfrag}
\usepackage{url}
\usepackage{amsmath,amssymb,xcolor}
\begin{document}
\title{Symmetry in  Self-Similarity in Space and Time---Short Time Transients 
and Power-Law Spatial~Asymptote}
\author{Ken Sekimoto}
\email{ken.sekimoto@espci.fr, corresponding author.}
\affiliation{Mati\`{e}res et Syst\`{e}mes Complexes, CNRS-UMR7057, Universit\'e   Paris-Diderot, 75205 Paris, France}
\affiliation{PCT, UMR CNRS 7083 Gulliver, ESPCI ParisTech PSL Research University, 10 rue Vauquelin, 75005 Paris, France}
\author{Takahiko Fujita}
\email{rankstatistics@gmail.com}
\affiliation{Faculty of Science and Engineering, Department of Industrial and Systems Engineering, Chuo University, 1-13-27, Kasuga, Bunkyo, Tokyo, Japan}
\begin{abstract}
The {self-similarity {in space and time (hereafter self-similarity)}}, either deterministic or statistical, 
 is characterized by similarity exponents and a function of scaled variable, called the scaling function. 
In the present paper, we address mainly the {self-similarity} in the limit of early stage, as~opposed to the latter one, and also consider the scaling functions that decay or grow algebraically, as~opposed to the rapidly decaying functions such as Gaussian or error function. 
In particular, in~the case of simple diffusion, our symmetry analysis shows 
a mathematical mechanism by which the rapidly decaying scaling functions are 
generated by other polynomial scaling functions. While~the former is adapted to the {self-similarity } in the late-stage processes, the latter is adapted to the early stages.  This  paper sheds some light on the internal structure of the family of  self-similarities generated by a simple diffusion equation.
 Then,   we present an example of {self-similarity } for the late stage whose scaling function has power-law tail, and also several cases of {self-similarity } for the early stages. These examples show the utility of {self-similarity } to a wider range of phenomena other than
the late stage behaviors with rapidly decaying scaling functions.
\end{abstract}
\maketitle
\section{Introduction}
When the evolution of state shows some scaling behavior involving both space and time without fixed characteristic length or time, it is said to have {self-similarity \null{in space and time}} (hereafter, ``{self-similarity}'' for short) \cite{Barenblatt-book}. 
\null{Self-similarity is a symmetry of the time-dependent fields which remain invariant under certain scale transformations of space and time.}
Similarity and dimensional methods in mechanics has long been known (\cite{SEDOV1959} and references therein).
Dynamical critical phenomena~\cite{Hohenberg-Halperin}
or the late-stage dynamics of the first-order transition~\cite{Gunton-Droz1983} are among the typical examples that have been studied a lot. Deterministic cases have also been actively studied (e.g., see the book~\cite{Barenblatt-book}). The~ {{self-similarity}} is represented by some scaled variables, scaling exponents, and~scaling functions. 
Among the simplest examples of {self-similarity} is the smearing Gaussian distribution governed by a diffusion equation, where the diffusing field $\phi$ evolves, for~example, in~the form $\phi(x,t)=t^{-\alpha}f(xt^{-\inv{2}}),$ where $f(s)=\exp(-\mbox{const.}\times s^2)$ is the scaling function for the scaled variable, $xt^{-\inv{2}},$ and $\alpha$ is the scaling exponent, being $1/2$ in this case.  \null{Such symmetry object, scaling function, can have different types, and~one of our main goals is to relate the scaling functions of different types.}

Thus far,   most studies on the  {{self-similarity}} in physics have   focused on the late-stage evolution, i.e.,~well after the initial transient, and~at the same time focused on the case where the scaling function or its derivative decays rapidly, i.e.,~faster than any power-law as a function of the scaled variable. 
The basic reason for this is the fact that the group velocity in the parabolic or hyperbolic evolution equation is vanishing or finite at large length scale. This implies that the  {{self-similarity}} with algebraic asymptotes of the scaling function is reduced to   {\it {static} 
} \null{(with steady flux)} 
algebraic tails in the far field.   As~such  far field condition has been rarely prepared physically, we have   not paid much attention to the 
 {{self-similarity}} with power-low asymptotes in the scaling~function.

  {{Self-similarity}} does not exclusively occur in the late stage, where the spatial scale is also large.     {{Self-similarity}} as the short-time limit has not drawn much attention thus far, to~the authors' knowledge. In~this limit, the~local region is infinitely magnified by the scaled variable. Then,   the~scaling functions are not constrained to decay rapidly but they can even grow algebraically for large values of scaled~variable.

In this paper, we first analyze the internal symmetry properties in the family of scaling functions for the 1D diffusion system {(Section \ref{sec2})}.  Through the equation governing the scaling function and its adjoint equation together with the operator of ``Wick rotation'',
 the scaling functions with rapidly decaying asymptote, which are suitable for the late-stage {self-similarity}, are related to the polynomial scaling functions, suitable for the short-time {self-similarity}. 
Secondly, we present some examples of {{self-similarity}} having the scaling functions with algebraic asymptote {(Section \ref{sec:algebtailed})}:   The first example is the late-stage {{self-similarity}} having an algebraically decaying scaling function which we think to be realizable in the permeative relaxation process of hydrogel. The~other examples are the short-time {{self-similarity}} of 1D diffusion from parabolic or cusp-like initial field configuration  and 
the short-time {{self-similarity}} of 1D linearized capillary-driven thin-film equation.
Through these studies, we shed light on the family of {{self-similarity}} behaviors and the relationship among them. In~  {Section \ref{sec4},} we~briefly discuss the generalizability of the  results obtained in this~paper.

\section{Symmetry Properties in the Family of Scaling Functions for the 1D Diffusion~Process}\label{sec2}
\vspace{-6pt}
\subsection{Scaling~Functions}
For the 1D diffusion, $\partial u/\partial t=\partial^2 u/\partial x^2,$ 
\null{different approaches have been developed (see, for~example,~\cite{group-anal-heat-eq1999}). Here, we focus on the solutions with self-similarity in space and time}. 
\eqn{\label{eq:defUp} 
u(x,t)=u_p(x,t)\equiv {t^{-{\,(p+1)}/{2}}} \phi_p({t^{-1/{2}}{x}}),}
where the scaling function $\phi_p(s)$ should obey
\eqn{\label{eq:opL}
\mathcal{L}_p\, \phi_p(s)=0, \mbox{with} \mathcal{L}_p=\frac{d^2}{ds^2}+\frac{s}{2}\frac{d}{ds}+\frac{p+1}{2}.
}   

The index $p$ has been adjusted so that $p=0$ corresponds to the Gaussian scaling function.
The even and odd solutions, $\Phi_{p}^{\rm (even)}(s)$ and $\Phi_{p}^{\rm (odd)}(s),$ respectively, are given in Appendix \ref{app:family}
 with some explicit form for several $p$ values. Here,   we note only those properties that are relevant for the discussion below:
For non-integer values of $p$, the~solution behaves asymptotically with power law, $|\phi_p(s)|\sim |s|^{-(p+1)},$ for $|s|\to \infty.$ The special cases are: (i) the solutions with rapidly decaying derivatives for non-negative integer $p$ with even [odd] parities when p is even [odd]; and (ii) the polynomial solutions for negative integer $p$ with even [odd] parities when p is odd [even]. Surprisingly, (i) and (ii) are related to each other, as~shown below. \null{For the diffusion equation or heat equation, the~general similarity solutions have been studied in detail (see, for~example,~\cite{bluman-cole-1969}). However, the~focus has not been put on the relationship among the solutions.}

\subsection{Symmetry in the Family of Scaling~Functions}
We noticed the two independent ways to relate the operators 
$\mathcal{L}_p$ defined by Equation~(\ref{eq:opL}) for $p\in \R$ and their adjoint ones,
$\mathcal{L}^\dag_p$, defined by
 \eqn{
 \mathcal{L}_p^\dag=\frac{d^2}{ds^2}-\frac{1}{2}\frac{d}{ds}s+\frac{p+1}{2}.}
The first relation between $\mathcal{L}_p$ and $\mathcal{L}^\dag_p$ is nothing but 
the definition of adjoint, i.e.,~the ``integration by parts'' identity,
\begin{widetext}
\beqa \label{eq:dual}
&&\psi_p(s) \mathcal{L}_p\phi_p(s)-\phi_p(s)\mathcal{L}_p^\dag \psi_p(s)
=\frac{d}{ds}\inSbracket{ \psi_p'(s) \phi_p(s)-\psi_p(s) \phi_p'(s)-\frac{s}{2}\psi_p(s) \phi_p(s)}.
\eeqa
\end{widetext}
 Equation~(\ref{eq:dual}) means that $\phi_p$ satisfying $\mathcal{L}_p\phi_p=0$ 
and $\psi_p$ satisfying $\mathcal{L}^\dag_p \psi_p(s)=0$ 
are mutually the integrating factor of the other. In~fact, if~we have such $\psi_p,$ 
the solution of $\mathcal{L}_p\phi_p=0$ can be constituted from the vanishing condition of the right-hand side of Equation~(\ref{eq:dual}).
The latter condition leads to 
\beq \label{eq:integrated}
\phi_p(s)= C_1\, \psi_p(s)e^{-\frac{s^2}{4}}+C_2 \int^s_{1}\frac{\psi_p(s)}{[\psi_p(u)]^2}e^{-\frac{s^2}{4}+\frac{u^2}{4}}du,
\eeq
where $C_1$ and $C_2$ are the integration constants, and~the integration pathway is assumed to be properly chosen to avoid the poles of the integrand.
The second relation between $\mathcal{L}_p$ and $\mathcal{L}^\dag_p$ is through a kind of Wick rotation:
\eqn{\label{eq:hidden}
 {\mathcal{L}}_{p}\circ J = - J\circ \mathcal{L}^\dag_{q}, \qquad p+q+1=0, 
 }
where the operator $J$ is defined such that $J\circ \Phi(s)=\Phi(is)$ for any object including $s$ as a variable. The~identity in Equation~(\ref{eq:hidden}) can be verified by   elementary~calculus.

The presence of two different pathways connecting the family  $\{\phi_p\}$ with that of 
$\{\psi_q\}$ reveals the intra-relationship among the scaling functions $\{\phi_p\}$, which is mediated by $\{\psi_q\}.$ Figure~\ref{fig:Wick} schematically shows this mechanism:
\begin{figure}[]\centering %
\includegraphics[width=1.3 in,angle=-90]{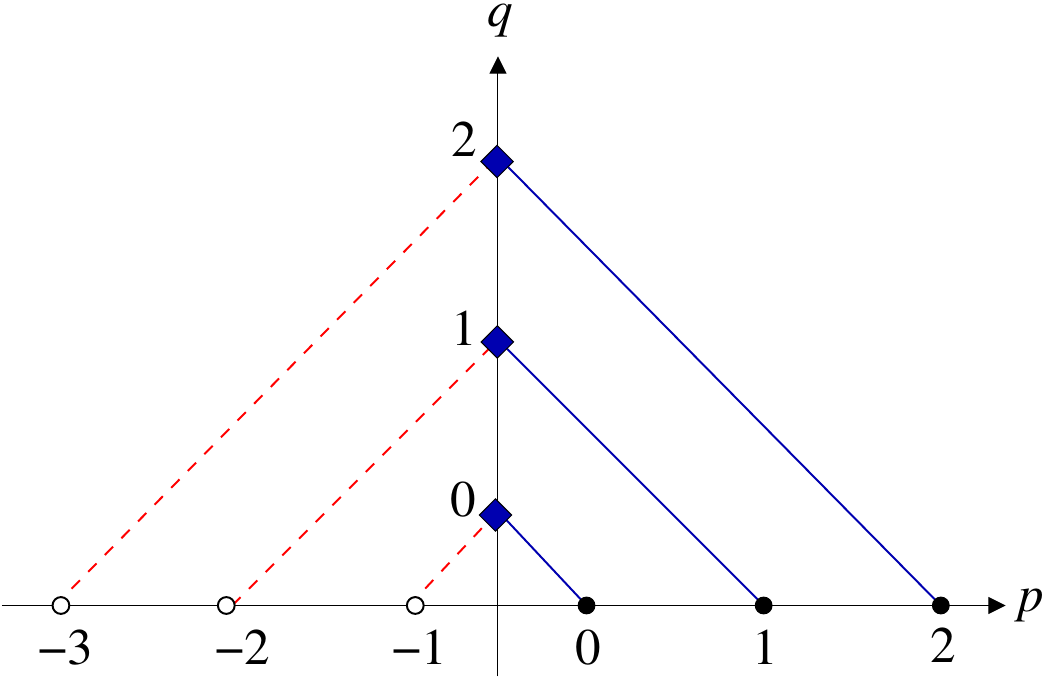}
\caption{\null  Diagram showing the relations among the solutions of $\mathcal{L}_p$ 
and those of $\mathcal{L}^\dag_q.$  See the main text.} 
\label{fig:Wick} 
\end{figure}
A solution of $\mathcal{L}_p\phi_p=0$ (e.g.,    $p=-3$) gives the solution
of $\mathcal{L}^\dag_q\psi_q=0$ 
 with $p+q+1=0$ (e.g.,    $q=+2$) through $\phi_{p}(s)=\psi_q(is).$
Once $\psi_q(s)$ satisfying   $\mathcal{L}^\dag_q\psi_q=0$ is given, 
the solution of $\mathcal{L}_{p}\phi_{p}(s)=0$ with $p=q$ (e.g.,    $p=2$)
is constructed by Equation~(\ref{eq:integrated}) with $\psi_p(s)$ and $\psi_p(u)$ being replaced by $\phi_{-(p+1)}(-is)$ and  $\phi_{-(p+1)}(-iu),$ respectively.

While the above symmetry property applies for any real index $p$, 
the cases of particular interest are when $\phi_{-p-1}(s)$ are polynomial with  non-negative integers $p.$ Such $\phi_{-p-1}(s)$
can generate, through~$\psi_p$, the~$\phi_{p}(s)$ whose derivative decays rapidly 
as the term $\propto C_1$ in Equation~(\ref{eq:integrated}).
Such polynomial solution $\phi_{-p-1}(s)$ has its origin in the adjoint equation $\mathcal{L}^\dag_p\psi_p=0,$ which can in turn be reduced to the Hermite equation by changing variable, $y=s/2.$ This type of analysis may shed some light on the {\it internal} structure of the family of dynamical similarities of a system of evolution, although~for the moment being limited to a particular case.
In Section~\ref{subsec:apexes}, we study the case where $\phi_{-p-1}$ is $\Phi_{-3,e}$ or 
$\Phi_{-2,o}$ as short-time self-similar evolution. The~above-mentioned symmetry connects
these scaling functions to the rapidly decaying ones $\Phi_{2,e}$ and $\Phi_{1,o}$ in the 
notation used in Appendix \ref{app:family}.

\section{\null{Self-Similarity in Space and~Time} with Scaling Function Having Algebraic Tail}
\label{sec:algebtailed}
As mentioned above, we should realize that, when the scaling function of {self-similarity} has algebraic tails as function of the scaled variable, these algebraic asymptotes must be realized in the {\it initial} field configuration. 
For example, in~Equation~(\ref{eq:opL}) for the 1D diffusion equation, the~algebraic tail issues from the last two terms in $\mathcal{L}_p,$ whose balance implies  $\partial u/\partial t\simeq 0.$ Indeed, if~we put  $\sim |s|^{-(p+1)}$ for $\phi_p(s)$ in Equation~(\ref{eq:defUp}), we have $u_p(x,t)\sim t^0\, x^{-(p+1)}.$
This consequence is intuitively understandable because {most of the evolution system of our interest spreads mainly locally, being unable  to update substantially the configurations far apart.}
Below, we present some examples of {self-similarity} having algebraic spatial tail.
One is for the late stage and the other two are for the short-time~limit.

\subsection{Late-Stage {Self-Similarity \null{in Space and~Time}} of Gel Network\label{subsec:gel}}
\begin{figure}[]\centering %
\includegraphics[width=2.0 in]{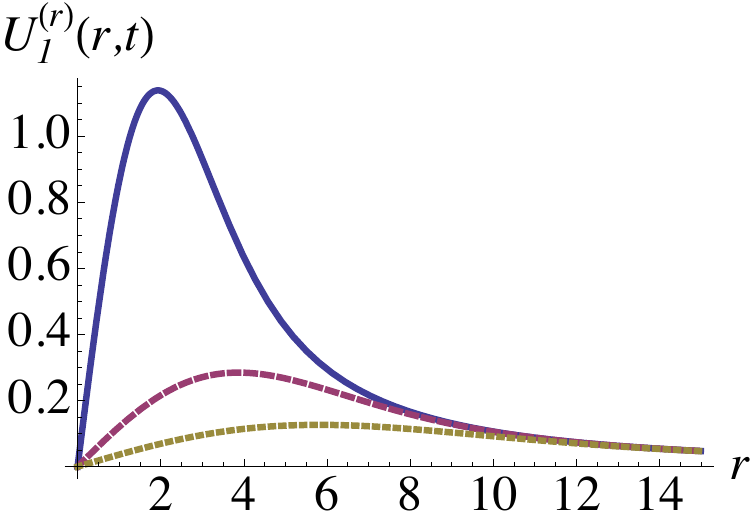}
\caption{  
\null {Radial displacement of the gel.}
The radial displacement, $U_{\rm rad}(r,t)=U_1^{(r)}(r,t)\equiv 
t^{-1}\Psi_1(r/\sqrt{t}),$ is shown  vs.  $r$ 
for the various time, $t=1$ (blue), $4$ (magenta), and~$9$ (yellow). 
The vertical scale has been arbitrarily chosen and $R_0$ is supposed to be infinitesimal.
  For $r\gg \sqrt{t}$, the~displacement, ${U_{1}^{(r)}}(r,t),$ retains a time-independent tail, $\sim r^{-2}.$ 
While the finiteness of $R_0$ should modify the radial displacements at small scale, $r\sim R_0,$ it does not influence the fate of the algebraic tail.} 
\label{fig:Ut0-t149} 
\end{figure}

It is generally difficult to prepare an algebraic static asymptotes over a large spatial range. This~may be the reason the {self-similarity} having algebraic tails has been scarcely studied in physical application. 
Below, we   present an example that an algebraic tail is nevertheless thinkable.
The system considered is an aqueous gel, i.e.,~a soft 
solid consisting of a very fine permeable network saturated with water.
Such a system responds as if it were an incompressible isotropic elastic body at short-time scale while its permeative relaxation under an osmotic pressure gradient is very slow. 
As the setup, we suppose a uniform 3D aqueous gel having a small water-filled hole around the center, $r=0.$ The~hole's radius is $R_0.$ 
For $t<0$, the~gel is in uniform equilibrium under no stress and, at~$t=0$, we quickly inject a small amount of water into the hole. The~injection makes expand the holes from radius $R_0$ to $R_1 (>R_0).$ 
Soon after the injection, a~quick elastic re-equilibration without permeation will take place, where each material point on the gel network has been pushed off radially. We denote by $U_{\rm rad}(r,0^+)$ the initial radial displacement of the material point which has been at the radius $r.$ This displacement must satisfy $4\pi r^2 U_{\rm rad}(r,0^+)=({4\pi}/3)[R_1^3-R_0^3]$ because the gel network behaves incompressible at short timescale (\null{where we assume  $|U_{\rm rad}|\ll R_0$}).
Thus, we can prepare an algebraic initial configuration, $U_{\rm rad}(r,0^+)\propto r^{-2},$ which is \null{exact in the limit of for $r/R_0\to \infty.$}
The subsequent very slow permeative relaxation obeys the standard elasto-hydrodynamic model~\cite{THB}, which gives the following time evolution equation for the radial displacement field $U_{\rm rad}(r,t)$:
\eqn{\label{eq:Deq3}
\inRbracket{\frac{\partial}{\partial t}-D
\inSbracket{\frac{\partial^2}{\partial r^2} -\frac{2}{r^2} } }[r \,U_{\rm rad}(r,t)]=0,
}
where $D$ is the (cooperative) diffusion constant $D=(K+4\mu/3)/f,$ with $f$, $\mu$, and~$K$ being, respectively, the~permeative friction coefficient between the gel network and solvent, the~shear modulus, and~the bulk (osmotic) modulus of the gel~\cite{deGennes}.
Hereafter, we   use the units such that $D=1.$
In the late-stage, $\sqrt{t}\gg R_0,$ at large distance $r\gg R_0,$ 
Equation~(\ref{eq:Deq3}) allows the similarity solutions of the form, $U_{\rm rad}(r,t)=U^{(r)}_{p}(r,t)\equiv {t^{-{(p+1)}/{2}}}\Psi_p\inRbracket{{r}/{\sqrt{t}}},$ where 
$p=1$ is compatible with  our initial condition. 
See {Appendix~\ref{sec:app-gel3D}} for the details including the analytical solutions. 
Figure~\ref{fig:Ut0-t149} shows the evolution of {self-similarity} solution, $U_{1}^{\rm (r)}(r,t)=t^{-1}\Psi_1(r/\sqrt{t}),$ at $t=1,4$ and $9.$ 
We see there that the algebraic tail of the radial displacement at large $r$ remains unchanged until the time $t\sim r^2.$
\unskip
\subsection{\null{Self-Similarity in Space and~Time} in the Short-Time Limit of 1D Diffusion}
\label{subsec:apexes}
\begin{figure}[]
\centering     
\subfigure[\null]{\label{fig:round-exact}  \includegraphics[width=4cm]{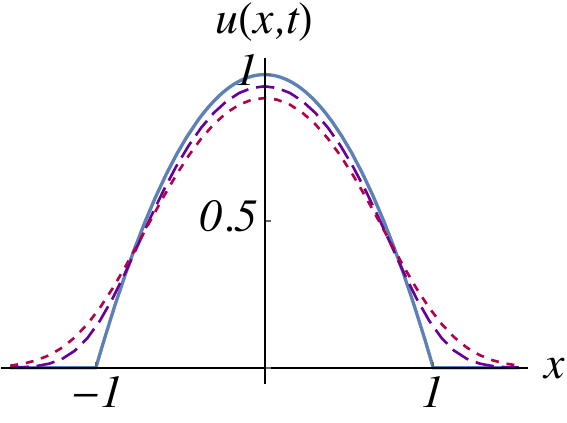}}
\hspace{.5cm}  
\subfigure[\null]{\label{fig:wedge-vs-x}  \includegraphics[width=4cm]{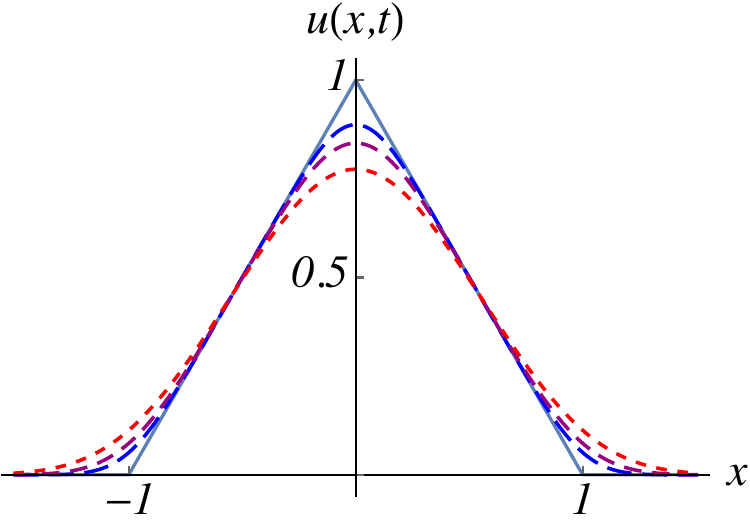}}
\hspace{3.5cm}    
\subfigure[\null]{\label{fig:round-fit}  \includegraphics[width=4cm]{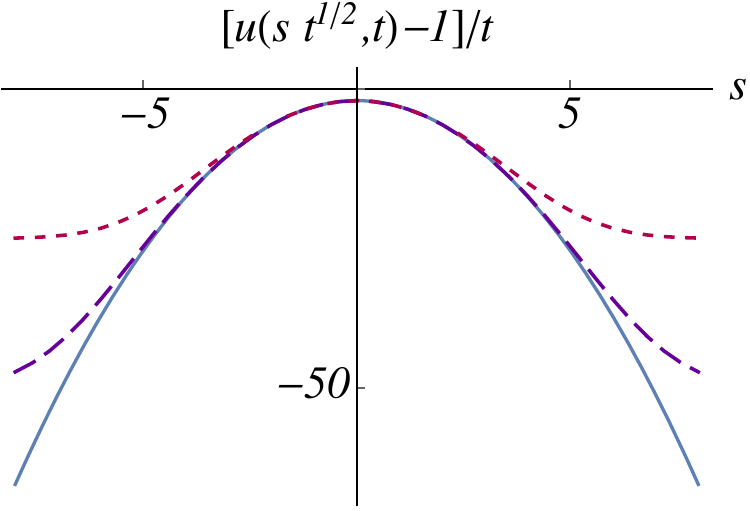}}
\hspace{0.5cm}  
\subfigure[\null]{\label{fig:cusp-fit}  \includegraphics[width=4cm]{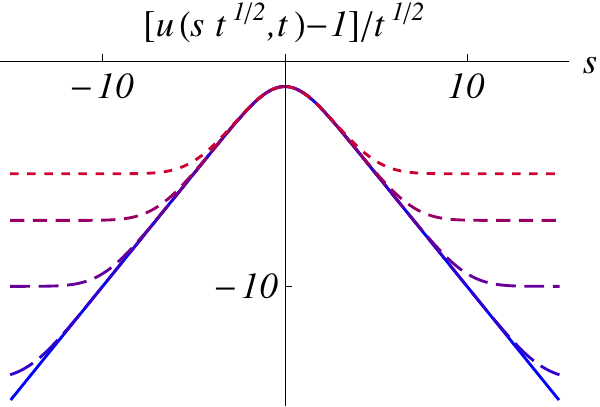}}
\caption{\null 
\label{fig:wedge-rounding}
(\textbf{a},\textbf{b})  Initial conditions (solid curves) as well as the short-time evolution at
$t=0.02$ and $t=0.04$ (dashed curves in (\textbf{a}))~
and at  $t=0.005, 0.01,$ and $0.02$ (dashed curves in (\textbf{b})),~
respectively, are shown versus the unscaled variables. 
The initial condition has $u(0,0)=1$ and the support, $[-1,1].$
(\textbf{c},\textbf{d})   The scaling function and the short-time evolution are represented using the scaled variables:
(\textbf{c}) replots the evolution in (\textbf{a}) using  $[u(x,t)-u(0,0)]/t$ and $s=x/\sqrt{t}$ and compares with $\phi_{-3}(s)=2+s^2,$ while (\textbf{d}) replots the evolution in (\textbf{b}) using $[u(x,t)-u(0,0)]/\sqrt{t}$ with $s$ and compares with $\phi_{-2}(s)=e^{-{s^2}/{4}}+({\sqrt{\pi}\,s}/{2}){\rm erf}({s}/{2}).$ 
} 
\label{fig:wedge}
\end{figure}

Parabolic or wedge-shaped apexes are among the most generic {\it local} configurations that we encounter as the initial conditions of diffusion. While the former is generic form around the smooth extrema, the~latter should be realized when some spatially localized source is removed after having established the wedge-shaped steady state profile. 
We discuss the simple case of one-dimension, but~the generalization to higher dimensions is straightforward. 
In {Figure} \ref{fig:wedge}a,b, the~parabolic and wedge-like initial conditions are shown by solid curves, respectively. Their functional forms are
 $u(x,0)=\max(1-x^2,0)$ and $u(x,0)={\max}[1-|x|,0]$. 
 The time evolution in the early stage is also shown in the same figure.
In Figure~\ref{fig:wedge}c,d, the~same evolutions are replotted by the scaled variables.
In both cases, parabolic or wedge-like apexes, we~can see that, the~earlier is the time, the~wider is the range of the scaled variable $s$ that fits the scaling~{functions.} 

As     discussed below, the~self-similarity in Figure~\ref{fig:wedge}a,c should be distinguished from the spreading of a Gaussian~peak.

\begin{figure}[]{}
\centering{
\includegraphics[width=1.5in]{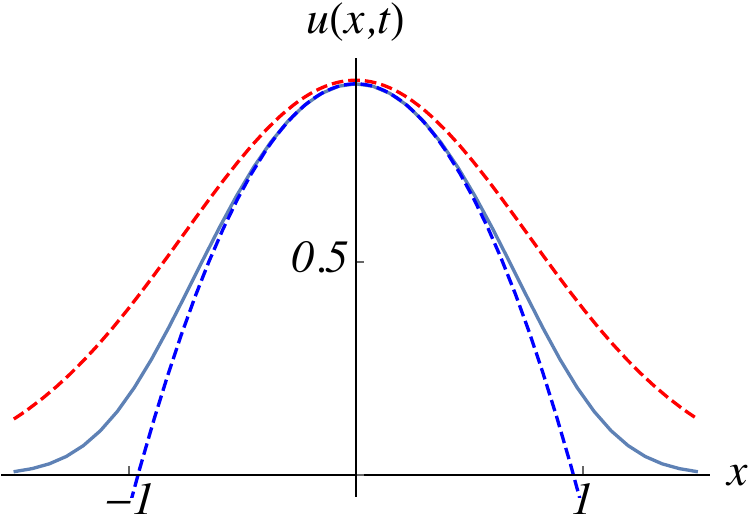}
}
\caption{\null The early stage profile at $t=0.04$ calculated from the truncated parabolic initial condition (solid curve) is compared with the short-time {self-similarity}  prediction (lower dashed curve) and with the late-stage one (upper dashed curve).
\label{fig:parabolla-shifting}} 
\end{figure}

The short-time {self-similarity} serves to determine the kinetic parameters in the evolution equations. In~the above example,  
the knowledge of the parabolic scaling function $\phi_3(s)=2+s^2$ allows
  fitting  the smooth apex of $u(x,t)$ in Figure~\ref{fig:wedge}a
with $u_{{3}}(x,t):=A[(x-B)^2+C+Dt],$ which then gives the diffusion coefficient $D$ 
among other parameters.
An important point is that the peak curvature $A$ remains constant in the short-time {self-similarity}, as~opposed to the late-stage {self-similarity}. In~the latter framework, the~Gaussian scaling function fitted to the initial parabolic apex predicts the decrement of the curvature linear in time. Figure~\ref{fig:parabolla-shifting} compares the 
two approaches.
 While the short-time {self-similarity} (lower dashed curve) keeps the initial curvature at the apex and fits well to the true evolution (solid curve), the~late-stage {self-similarity} (upper dashed curve) fails to fit the apex curvature. After~preparing the first version of our paper~\cite{SF-arxiv2012},
we noticed \citet{stonepnas} published at the same period, where the authors addressed the issue of early-stage  vs.~late-stage. 

\subsection{Short-Time {Self-Similarity \null{in Space and~Time}} of Capillary-Driven Thin-Film Equation}
The capillary evolution of the profile of the free surface of a thin viscous liquid film occurs as the quasi-static balance between the Laplace pressure and the viscous force due to the flow (see~\cite{Benza-epje2013}
 and the references cited therein).
The time evolution of the height profile function, $u(x,t),$ obeys the fourth-order partial differential equation:
\eqn{\label{eq:pde4}
\frac{\partial u}{\partial t}=-\frac{\partial^4 u}{\partial x^4} .}

The system allows the {self-similarity}, ${u(x,t)=t^{-\frac{\beta+1}{4}}\,\phi({x}{t^{-\inv{4}}}), }$ where the scaling function $\phi(s)$ with $s=x t^{-\inv{4}},$ should obey
\eqn{ 
-\phi''''(s) + \frac{s}{4}\phi'(s)  + \frac{\beta+1}{4}\phi(s) =0.
\label{eq:ODE4}}

\begin{figure}[]\centering %
\includegraphics[width=2.0 in]{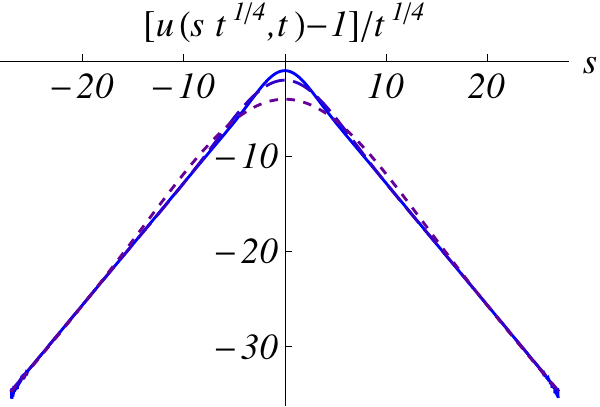}
\caption{\null {Short-time {self-similarity} evolution with algebraic asymptote of capillary-driven thin-film~equation.}
From the initial profile $ u(x,0)=1-|x|$ (this cannot be shown on this plane), the~solution $ u(x,t)$ for Equation~(\ref{eq:pde4}) 
at $t=1$ (solid curve),  $t=2^4$ (long dashed curve) and  $t=2^8$ (short dashed curve) are presented using the scaled variables. The~scaling function, which is numerically obtained by solving Equation~(\ref{eq:ODE4}) with $\beta=-2,$ is almost indistinguishable with the earliest profile (solid curve).
\label{fig:b1P0sur4} }
\end{figure}

Most generally the solution $\phi(s)$ or its derivative $\phi'(s)$ do not decay rapidly for a given value of $\beta.$ 
 Equation~(\ref{eq:ODE4}) under a given $\beta$ allows two independent even solutions of which one solution can have the power-low asymptote $\phi(s)\sim |s|^{-(\beta+1) }.$ In~\cite{Benza-epje2013}, the~special case corresponding to $\beta=0$ with rapidly decaying derivative has been analyzed. 
As in the case of diffusion, the~scaling function with power-law asymptote implies the static algebraic asymptote in space, $u(x,t)\sim |x|^{-(\beta+1)}$ ($|x|\gg t$), because~
 the algebraic asymptote of $\phi(s)$ issues from the balance of the last two terms on the left-hand side of Equation~(\ref{eq:ODE4}), which in turn implies $\partial u/\partial t\simeq 0$.
Figure~\ref{fig:b1P0sur4} shows, in~the scaled coordinates, the~evolution from the wedge-like initial configuration, which corresponds to $\beta=-2.$   
Unlike the simple diffusion, identifying the solution with algebraic asymptote 
requires a tuning of the ratio $\phi''(0)/\phi(0)$ (e.g.,    at $0.739\ldots$ for Figure~\ref{fig:b1P0sur4}) because the other solution of the same parity grows with $|s|$ faster than any power-laws. 
Zhang and Lister~\cite{similarity-rupture1999}, in~their study of the rupture of thin film by the  van der Waals interaction, discussed how the inclusion of nonlinearity in the evolution equation selects the value of the exponent $\beta$ to some discrete ones.

\section{Conclusions} 
\label{sec4}
{The fact that solutions to many differential equations enjoy the property of self-similarity is well known for about hundred years (see, e.g.,~the book by {L.I. }Sedov~\cite{SEDOV1959}).}
\null{In this paper, we   show  the existence of transformation that internally relates among the solution family of the scaling function (Figure  \ref{fig:Wick}). 
The concrete transformation consists of an integrating factor, which is the solution of adjoint equation and a Wick-rotation type transformation.}

To the authors' knowledge, this type of symmetry has not been reported before.
 How general is such structure for the evolution systems other than the diffusion remains a future problem. \null{The~existence of such transformation might not be specific for the differential equations of scaling functions. However, it is physically interesting and new that this transformation relates the self-similar relaxation process with flux-free boundary  on the one hand (i.e.,   $\phi_p$ with $p\ge 0$) and the self-similar process with steady flux at boundary on the other hand ($\phi_p$ with $p\le -1$).}

We   show  the utility of the \null{self-similarity \null{in space and time}}  to analyze the early stage of diffusive processes including the pure diffusion as well as the thin-film equation. These types of self-similarity focus on the local process in space-time and, therefore, are distinguished from the late-stage ones.
The merit of the short-time similarity is that, {even if the whole system has characteristic lengths or times, the~self-similarity in the local space-time region works at the scale which is far below those characteristic scales .}
\null{In Section~\ref{sec:algebtailed}, we   show  the utility of the scaling functions that decay or grow algebraically. The~short-time self-similarity of diffusion from smooth peak or wedge-like apex are new physical results, and~so is the short-time self-similar solution of the surface diffusion. Unlike~that of \citet{Benza-epje2013}, the~latter solution has algebraic tails.
Generically, one can prepare any algebraically decaying tail  as initial condition for which the approach taken in Section~\ref{subsec:apexes} will be useful, or~at least usable.}
\null{In the application of self-similar solution with algebraic tail to the relaxation of gel (Section \ref{subsec:gel}), we  bring a new physical idea of combining 
the elastic aspect of gel, which is described by hyperbolic or elliptic equations, and~the diffusive aspect, which obeys parabolic equations. By~the former aspect, we can prepare an algebraic tail, and~afterwards the gel shows self-similar relaxation having algebraic tail.}

{ In some fields of research~\cite{SEDOV1959},
 the main idea of considering scaling relations is the use of these relations for constructing unknown solutions to the equations whose physical behavior is yet to be clarified. In~the present paper, we do  the opposite. We  take  simple equations with well or less known solutions and state that, for~some specially constructed initial conditions,    scaling with respect to space and time variables occurs. In~the  future, these two approaches might meet, where the internal relationship among the solution family will help to understand the global picture.}

Lastly, the~application might not be limited to the linear systems. One apparent example is the Burgers equation, since this equation can be reduced to the diffusion equation through the Cole--Hopf transformation if the boundary conditions are compatible~\cite{cole-hopf}. Other type of generalization remains as an open question. \null{When some mathematical objects represents well some physical phenomena,
 the~mathematical relationship among the former might well have some physical consequences in the latter. We would like to know in the future the meaning of the transformation in Section \ref{sec2} in physics~terms.} 

\mbox{}

\appendix
\section{Family of Scaling Functions for the {Self-Similarity} Solution of 1D Diffusion~Equation}\label{app:family} 

{{When} 
 $q$ is non-negative integer, 
the adjoint equation $\mathcal{L}^\dag_q\psi_q(s)=0$
has a polynomial solution,} $\psi_q(s)=H_q({s}/{2})$ up to a multiplicative constant, where $H_q(x)$ is the $q$th-order Hermite {polynomial} function. By~the symmetry presented in the main text, the~polynomial solution of $\mathcal{L}_p\phi_p(s)=0$ with $p+q+1=0$ is  $\phi_{p}(s)=H_{q}({is}/{2}).$  
{For non-integer or negative values of $q,$
the symmetric solution satisfying ${\phi}_p^{\rm (s)}(0)=1,$ and ${\phi}_p^{\rm (s)\prime}(0)=0$ is
$ {\phi}_p^{\rm (s)}
(s)=e^{-\frac{s^2}{4}} {_1\! F_1}\left(-\frac{p}{2},\inv{2},\frac{s^2}{4}\right) $
and the antisymmetric solution satisfying ${\phi}_p^{\rm (a)}(0)=0,$ and ${\phi}_p^{ \rm (a)\prime}(0)=1$ is
${\phi}_p^{\rm (a)}(s)= s\,e^{-\frac{s^2}{4}} {_1\! F_1}\left(\frac{1-p}{2};\frac{3}{2};\frac{s^2}{4}\right)
,$
where  ${_1\! F_1}(a;b;z)$ is the confluent hypergeometric function of the first kind~\cite{conflhypergeom-Math,Abramowitz-Stegun}. }

\begin{widetext}
\begin{table}[]
\caption{}
\centering
\begin{tabular}{cll}
\toprule
\boldmath{${\emph{\textbf{p}}} 
$}	& \boldmath{$\Phi_{p,e}(s)$}	& \boldmath{$\Phi_{p,o}(s)$}\\
\hline
$-$3		& $\Phi_{-3,e}(s)=\frac{s^2}{2}+1$	& 
$\Phi_{-3,o}(s)= \frac{\sqrt{\pi}}{2}(\frac{s^2}{2}+1)\,{\rm erf}(\frac{s}{2})+\frac{s}{2}e^{-\frac{s^2}{4}}$
\\
$-$2		& $ {\Phi_{-2,e}(s)= \frac{\sqrt{\pi}}{2}
s\,{\rm erf}(\frac{s}{2})+e^{-{s^2}/{4}} }$ & $\Phi_{-2,o}(s)=s $\\
$-$1		& $\Phi_{-1,e}(s)=1$	& $\Phi_{-1,o}=\sqrt{\pi}{\rm erf}(\frac{s}{2})$\\
0		& $\Phi_{0,e}(s)=e^{-{s^2}/{4}}$\, (*) & $\Phi_{0,o}(s)= e^{-{s^2}/{4}}\int_0^s e^{{w^2}/{4}}dw$\\
1		& $\Phi_{1,e}(s)\propto\Phi'_{0,o}(s)$ 			& $\Phi_{1,o}(s)\propto\Phi'_{0,e}(s)$\, (*) \\
2		& $\Phi_{2,e}(s)\propto\Phi''_{0,e}(s)$\, (*)			& $\Phi_{2,o}(s)\propto\Phi''_{0,o}(s)$\\
\hline
\end{tabular}
\end{table}
\end{widetext}
\unskip
\begin{figure}[]\centering %
\includegraphics[width=1.8 in]{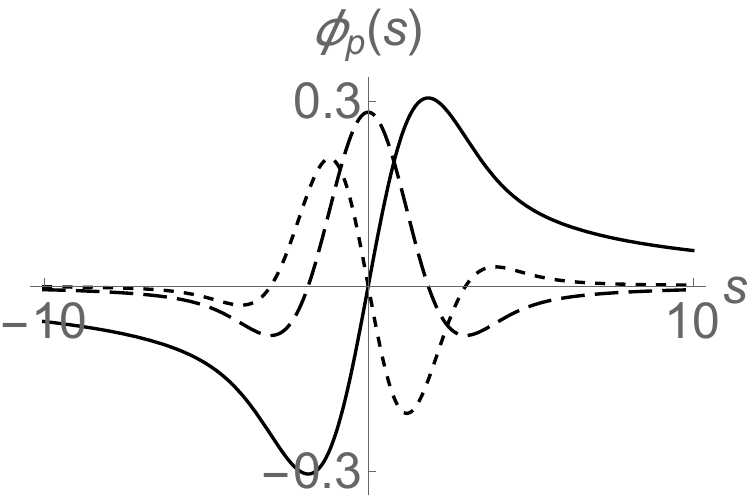}
\caption{\null {{Examples of} 
 the scaling functions satisfying Equation~(\ref{eq:opL}) with \null{power-law asymptotes}.} 
$\phi_{0,o}(s)= e^{-{s^2}/{4}}\int_0^s e^{{w^2}/{4}}dw$ (solid curve), 
$\phi_{1,e}(s)=\phi_{0,o}'(s)$ (coarse dashed curve), and~$\phi_{2,o}(s)=\phi_{0,o}''(s)$ (fine dashed curve). For~$|s|\gg 1$, these functions decays as $\sim |s|^{-(p+1)}$ with $p=0,1,$ and $2$, respectively.
}
\label{fig:phit-p012} 
\end{figure}

\section{Spacetime Self-Similarity Solution of Equation \eqref{eq:Deq3}  \label{sec:app-gel3D}}
Into Equation \eqref{eq:Deq3}, we substitute the form of the radial displacement, 
\[U_{}^{\rm (r,t)}=U_p^{\rm (r)}(r,t)\equiv\inv{(Dt)^{\frac{p+1}{2}}}\Psi_p\inRbracket{\frac{r}{\sqrt{Dt}}}.\] 

{We then }obtain the ordinary differential equation for $\Psi_p(s)$ with $s> 0$:
\eqn{\label{eq:psip}
2s^2 \Psi_p''(s)+s(s^2+4) \Psi_p'(s) +((p+1)s^2-4) \Psi_p(s)=0 .}

In addition to the far field condition,  $\Psi_p(\infty)=0,$ we impose the boundary condition 
at $s=0$, which is compatible with the homogeneous dilatation around the origin. 
In terms of the displacement vector field  $\vec{U},$ the homogeneous dilatation reads $\vec{U}(\vec{r},t) =c(t) \, \vec{r} +{\cal O}(r^3)$ for $r\downarrow 0$ with $c(t)$ being a function of time. We then impose $\Psi_p(0)=0.$ 
Under these boundary conditions, the~solutions of Equation \eqref{eq:psip} is written in terms of the  confluent hypergeometric function $_1\! F_1(a;b;z)$ \cite{conflhypergeom-Math,Abramowitz-Stegun}. 
The~solution reads
 \eqn{
\Psi_p(s)=\Psi'_p(0)\, s \times\, _1\! F_1\inRbracket{1+\frac{p}{2};\frac{5}{2};-\frac{s^2}{2}}.}
 
 This solution behaves as $\Psi_p(s)\simeq \Psi'_p(0)\times s $ for $0\le s\ll 1$ and $\Psi_p(s)\sim s^{-(p+1)}$ for $s\gg 1.$ Because~$U_{\rm rad}(r, 0)\sim r^{-2}$ for $r\gg t^{\inv{2}},$ we choose  $p=1,$  then we have
\eqn{\Psi_1(s)= \Psi'_1(0)\times \frac{6}{s^2}\inSbracket{-s\,e^{-\frac{s^2}{4}}+\sqrt{\pi}{\rm erf}\inRbracket{\frac{s}{2}}}
\label{eq:drho}.}



\end{document}